\documentclass[12pt]{amsart}
\usepackage{latexsym,amssymb,amsxtra,amsmath,amsfonts}
\usepackage{verbatim}
\usepackage{epsfig}
\usepackage[pdftex]{color} 
\usepackage[top=1in, bottom=1in, left=1in, right=1in]{geometry}

\newtheorem{theorem}{Theorem}
\newtheorem{lemma}[theorem]{Lemma}
\newtheorem{proposition}[theorem]{Proposition}
\newtheorem{corollary}[theorem]{Corollary}

\theoremstyle{definition}

\newcommand{\beq}{\begin{equation}}
\newcommand{\eeq}{\end{equation}}
\newcommand{\beqa}{\begin{eqnarray}}
\newcommand{\eeqa}{\end{eqnarray}}
\newcommand{\beaa}{\begin{eqnarray*}}
\newcommand{\ben}{\begin{eqnarray*}}
\newcommand{\eaa}{\end{eqnarray*}}
\newcommand{\een}{\end{eqnarray*}}

\numberwithin{equation}{section}
\numberwithin{theorem}{section}

\begin{document}

\title[The phase factors in singularity theory ]
{The phase factors in singularity theory}
\author{Todor Milanov}
\address{Kavli IPMU (WPI) \\ The University of Tokyo \\ Kashiwa \\ Chiba 277-8583 \\ Japan}
\email{todor.milanov@ipmu.jp}
\thanks{{\em 2000 Math. Subj. Class.} 14D05, 14N35, 17B69}

\begin{abstract}
The paper \cite{BM} proposed a construction of a twisted
representation of the lattice vertex algebra corresponding to the
Milnor lattice of a simple singularity. The main difficulty in
extending the above construction to an arbitrary isolated singularity is in the
so called {\em phase factors} -- the scalar functions produced by 
composing two vertex operators. They are certain family of multivalued
analytic functions on the space of miniversal deformations.  
The first result in this paper is an explicit formula for the unperturbed
phase factors in terms of the classical monodromy operator and the
polylogorithm functions. Our
second result is that with respect to the deformation parameters the
phase factors are analytic functions on the monodromy covering space. 

\end{abstract}

\maketitle
\tableofcontents
\addtocontents{toc}{\protect\setcounter{tocdepth}{1}}

\section{Introduction}
Lattice Vertex Operator Algebras (VOAs) were introduced by R. Borcherds in his work on
the Moonshine conjecture (see \cite{B1,B2}). After Borcherds' work, 
lattice VOAs were used as a tool to generalize the notion of a Kac-Moody
Lie algebra. In particular, it is very natural to apply Borcherds' construction to
singularity theory, because we have all the necessary ingredients: a Milnor
lattice equipped with the intersection pairing. What is even more
exciting is that there is a natural candidate for a root system,
namely the set of vanishing cycles. It is very tempting to investigate
to what extend the classical theory of simple Lie algebras extends to
singularity theory. We refer to \cite{Sa3, Sl} for further details on 
Lie algebras in singularity theory. The applications
of this class of Lie algebras is still not quite clear. In particular,
fundamental questions, such as finding the root multiplicities and
classifying the irreducible representations are wide
open. Nevertheless, Gromov--Witten theory provides a new
motivation to develop further Borcherds' ideas in the settings of
singularity theory. Especially, Givental's higher genus
reconstruction formalism \cite{G1,G3, Te} and the integrable
hierarchies introduced by 
Dubrovin and Zhang \cite{DZ} suggest to investigate whether the Lie
algebras in singularity theory are a source of integrability in the
same way as the ADE affine Lie algebras are a source of
integrability for the Kac--Wakimoto hierarchies \cite{KW}.

The relevance of period integrals in singularity theory for describing
integrable hierarchies was first observed by Givental in \cite{G3} and
developed further in \cite{FGM}. Partially motivated by these works,
it was suggested in \cite{BM} that the period integrals associated
with a given singularity can be used to construct a $\sigma$-twisted
representation of the lattice VOA associated with the
Milnor lattice, where $\sigma$ is the classical monodromy. The main
issue in confirming the prediction of \cite{BM} arises
when we compose two vertex operators. We get a set of phase factors
that we can not study only in the category of formal Laurent series, due
to the fact that exponentiating a formal Laurent series does not make
sense in general. The goal of this paper is to establish the analytic properties of the
phase factors that are necessary in order to extend the construction
of \cite{BM} in general. 

Our first result is that the vertex operators are mutually
local in the sense of VOA. The proof is based on an explicit
computation of the unperturbed phase factors. Our main tool is  a
formula due to Hertling that expresses Saito's higher residue
pairing in terms of the {\em Seifert  form}. The final answer is given in terms of the  polylogorithm
functions. Furthermore, it is not hard to see that the phase factors
are convergent and we can express them as integrals along the path of the 
so called {\em phase form} (see \eqref{pf}) -- quadratic form on the
vanishing homology with values 1-forms on the space of deformation
parameters. Using the locality and the fact that the
Gauss--Manin connection has regular singularities, we prove that the
phase factors are analytic functions on the monodromy covering
space. This fact was very important in both \cite{BM} and \cite{G3}
for constructing W-constraints  and Hirota quadratic equations of the
total descendant potential. The compatibility of the phase factors
with the monodromy representation is a very non-trivial and somewhat
surprising result.  It is equivalent to the following question of
Givental \cite{G3}: given a path $C$ and a pair of integral cycles
invariant along $C$, then is it true that the integral along $C$ of
the corresponding phase form is an integer multiple of $2\pi\sqrt{-1}$? 
Givental answered the above question positively in the case of simple
singularities. Unfortunately, his argument is hard to generalize,
because it relies on understanding the fundamental group of the
complement of the discriminant.

\section{Statement of the main results}

\subsection{Singularity theory.}
Let $W:(\mathbb{C}^{n+1},0)\to (\mathbb{C},0)$ be the germ of a
holomorphic function with an isolated singularity at $x=0$. We denote
by $x=(x_0,x_1,\dots,x_n)$ the standard coordinate system on
$\mathbb{C}^{n+1}$ and by 
\ben
H:=\mathcal{O}_{\mathbb{C}^{n+1},0}/(W_{x_0},W_{x_1},\dots,W_{x_n}),\quad
W_{x_i}=\partial W/\partial x_i,
\een
the Jacobi algebra of $W$. Let $\{\phi_j(x)\}_{j=1}^N$ be a set of
holomorphic germs that represent a basis of $H$. We may assume that
$\phi_1=1$. We are interested in the following deformation of $W$ 
\ben
F(x,t)=W(x)+\sum_{j=1}^N t_j\phi_j(x),\quad t=(t_1,\dots,t_N)\in \mathbb{C}^N.
\een
Let us cut a Stein domain $X\subset \mathbb{C}^{n+1}\times \mathbb{C}^N$ around
$(0,0)$, s.t., the map 
\ben
\varphi: X \to B\times \mathbb{C}_\epsilon,\quad (x,t)\mapsto (t, F(x,t) )
\een
is well defined and surjective, where $B\subset \mathbb{C}^N$ is a sufficiently small ball
around $0$ and $\mathbb{C}_\epsilon=\{\lambda\in \mathbb{C}:
|\lambda|<\epsilon\}$, and the fibers $X_{t,\lambda}:=\varphi^{-1}(t,\lambda)$ 
satisfy an appropriate transversality condition (see
\cite{AGV}). Put $S:=B\times \mathbb{C}_\epsilon $ to avoid cumbersome
notation.

Let $C$ be the critical locus $\{ F_{x_0}=\cdots = F_{x_n}=0\}\subset X$
of $F$ relative to $B$. The map
\ben
\mathcal{T}_B\to
\operatorname{pr}_*\mathcal{O}_C,\quad \partial/\partial
t_i\mapsto \partial F/\partial t_i \ (\operatorname{mod}\ F_{x_0},\dots,F_{x_n})
\een
is called the {\em Kodaira--Spencer} map. Decreasing $X,B$, and
$\epsilon$ if necessary we may arrange that the  Kodaira--Spencer map
is an isomorphism. In particular, each tangent space $T_tB$ has an
associative algebra structure. Let us fix a {\em primitive} form $\omega\in
\Omega^{n+1}_{X/B}(X)$, so that $B$ becomes a Frobenius manifold (see
\cite{He, SaT}). The flat Frobenius pairing is given by
the following residue pairing:  
\beq\label{res-pairing}
(\phi_1(x),\phi_2(x))_t := \frac{1}{(2\pi\sqrt{-1})^{n+1}} \int_{\Gamma}
\frac{\phi_1(y)\phi_2(y)}{F_{y_0}(t,y)\cdots F_{y_n}(t,x)} dy_0\dots dy_N,
\eeq
where the cycle $\Gamma$ is a disjoint union of sufficiently small
tori around the critical points of $F$ defined by equations of the type
$|F_{x_0}|=\cdots =|F_{x_n}|=\epsilon$ and $y=(y_0,\dots,y_n)$ is a
local coordinate system near each critical point, s.t.,
$\omega=dy_0\wedge \cdots \wedge dy_n$.  We have the following identifications:
\beq\label{flat-triv}
T^*B\cong TB\cong B\times T_0B\cong B\times H,
\eeq
where the first isomorphism is given by the residue pairing, the
second by the Levi--Civita connection of the flat residue pairing, and the last one is the
Kodaira--Spencer isomorphism. Let us denote by $\partial_i$ the flat
vector fields  corresponding to the basis $\{\phi_i\}_{i=1}^N\subset
H.$ 

\subsection{Period vectors}

Removing the singular fibers from $X$ we obtain a smooth 
fibration $X'\to S'$ known as the Milnor
fibration, where $S'$ is the subset
parametrizing non-singular fibers. Its complement is an irreducible
analytic hypersurface known as the {\em discriminant}. The homology $H_n(X_{t,\lambda};\mathbb{C})$
and the cohomology groups $H^n(X_{t,\lambda};\mathbb{C})$ form vector
bundles on $S'$ equipped with a flat Gauss--Manin connection. Let us fix a
reference point  $(0,\epsilon_0)\in S'$ and denote by 
\ben
\mathfrak{h}:=H^n(X_{0,\epsilon_0};\mathbb{C}),\quad
\mathfrak{h}^*:=H_n(X_{0,\epsilon_0};\mathbb{C})
\een
the reference fibers. The vector space $\mathfrak{h}^*$ contains the
so called {\em Milnor lattice} $Q:=H_n(X_{0,\epsilon_0};\mathbb{Z})$
and if we assume that $n$ is even then the intersection pairing gives
a symmetric bi-linear pairing on $Q$ (which however might be
degenerate). We normalize the intersection pairing by the sign
$(-1)^{n/2}$ in order to obtain a pairing $(\, |\, )$, s.t.,
$(\alpha|\alpha)=2$ for every {\em vanishing cycle} $\alpha\in Q$. We
refer again to \cite{AGV} for some more details on the vanishing
homology. 

Using the parallel transport with respect to the Gauss--Manin
connection we get the so called monodromy representation 
\beq\label{mon-rep}
\pi_1(S')\to \operatorname{O}(\mathfrak{h}^*),
\eeq
where the RHS denotes the group of linear transformations that preserve
the intersection pairing. The image of the monodromy representation is
a reflection group $W$ generated by the reflections
\ben
s_\alpha(x) = x-(\alpha|x)\alpha,
\een
where $\alpha$ is a vanishing cycle and $s_\alpha$ is the monodromy
transformation representing a simple loop corresponding to a path from
the reference point to a generic point on the discriminant along which
$\alpha$ vanishes. 

Let us introduce the notation $d_x$,
where $x=(x_1,\dots,x_m)$ is a coordinate system on some manifold,
for the de Rham differential in the coordinates $x$.
This notation is especially useful when we have to apply $d_x$ to
functions that might depend on other variables as well. 
The main object of our interest are the following period integrals
(see \cite{G3})
\beq\label{period}
I^{(k)}_\alpha(t,\lambda) = - d_t\ (2\pi)^{-1}\, \partial_\lambda^{k+1} \
\int_{\alpha_{t,\lambda}} d_x^{-1}\omega\ \in T_t^*B\cong H,
\eeq
where $\alpha\in \mathfrak{h}^*$ is a cycle from the vanishing homology,
$\alpha_{t,\lambda}\in H_n(X_{t,\lambda},\mathbb{C})$ is the parallel
transport of $\alpha$ along a reference path, 
and $d_x^{-1}\omega$ is a holomorphic $n$-form $\eta\in
\Omega^n_{X/B}$ defined in a neighborhood of the fiber
$X_{t,\lambda}$, s.t., $d_x\eta=\omega$. The periods are multivalued analytic
functions in $(t,\lambda)\in S$  with poles along the discriminant.
In other words they are analytic functions $\widetilde{S'}\to H$,
where $\widetilde{S'}$ is the monodromy covering space of $S'$. 

The periods satisfy the following system of differential equations
\begin{align}
\label{PF-1}
\partial_i I_\alpha^{(k)}(t,\lambda) & = -\phi_i \bullet_t I_\alpha^{(k+1)}(t,\lambda) \\
\label{PF-2}
\partial_\lambda I_\alpha^{(k)}(t,\lambda)  & = I_\alpha^{(k+1)}(t,\lambda) \\
\label{de-prim}
(\lambda-E\bullet_t)\partial_\lambda I_\alpha^{(k)}(t,\lambda)   & = 
\Big(\theta -k-\frac{1}{2}\Big) I_\alpha^{(k)}(t,\lambda),
\end{align} 
where the notation in \eqref{de-prim} is as follows. The vector field
$E\in \mathcal{T}_B$  is the {\em Euler} vector field, which by definition
corresponds to $F$ via the Kodaira spencer isomorphism. Changing the
basis $\{\phi_i\}_{i=1}^N\subset H$ if necessary we may arrange that
the Euler vector field takes the form
\beq\label{Euler-vf}
E=\sum_{i=1}^N (1-d_i)\tau_i\partial_i + \rho,
\eeq
where $\rho$ is a flat vector field and the degree spectrum
$1=d_1<d_2\leq \cdots \leq d_{N-1}<d_N=:D$ satisfies
$d_i+d_{N+1-i}=D$. The number $D$ is the {\em conformal dimension} of
the Frobenius structure. The linear operator $\theta$ is defined by 
\ben
\theta:H\to H,\quad \phi_i\mapsto (D/2-d_i)\phi_i.
\een
It is sometimes called the {\em Hodge-grading operator}.

\subsection{Extending the domain of the period vectors}
Using the differential equations \eqref{PF-1}--\eqref{de-prim} we
would like to extend the domain $S$ in such a way that the period vectors have a
translation symmetry and that there exists a point $u^{(0)}\in B$, such that
the singularities of the period vectors $I^{(m)}_\alpha(u^{(0)},\lambda)$
are at points $\lambda=u_i$ $(1\leq i\leq N)$ that coincide with the
vertices of a regular $N$-gon with center $0$. 

The differential equation \eqref{de-prim} allows us to extend the
definition of $I^{(m)}_\alpha(t,\lambda)$ for all $(t,\lambda)\in
B\times \mathbb{C}$. 
Let us denote by $t=(t_1,\dots,t_N)$ the standard coordinate system on
$B\subset\mathbb{C}^N$ induced from the linear coordinates on
$\mathbb{C}^N$. We denote by $t-\lambda\mathbf{1} =
(t_1-\lambda,t_2,\dots,t_N)$, i.e., this is the time $-\lambda$-flow
of $t$ with respect to the flat identity. Note that
$X_{t,\lambda}=X_{t-\lambda\mathbf{1},0}$, so the period vectors have
the following  translation symmetry
\beq\label{tr-sym}
I^{(k)}_\alpha(t,\lambda) = I^{(k)}_\alpha(t-\lambda\mathbf{1},0),
\eeq
where $|\lambda|\ll 1$. Note that the RHS is a multivalued function on
$B':=B\times\{0\}\cap S'$, while the LHS is a multivalued function on
$S'$. We choose $-\epsilon_0\mathbf{1}$ as a 
reference point in $B'$ and the reference path for the RHS is obtained
from the reference path for the LHS via the translation map
$$
\tau: S'\to B', \quad(t,\lambda)\mapsto t-\lambda\mathbf{1}.
$$ 
Using \eqref{tr-sym}, we extend the definition of the periods also for all $t$ that belong to 
a domain in $\mathbb{C}^N$ obtained from $B$ by  the flow of the flat vector field
$\partial_1$. Slightly abusing the notation we denote by $B$ the
extended domain and redefine $S$ to be $B\times \mathbb{C}.$ Note that
the Frobenius structure also extends on $B$ and now the translation
symmetry \eqref{tr-sym} makes sense for all $(t,\lambda)\in S$. 

Let us fix a point $t_0\in B$, s.t., the critical values of $F(x,t_0)$
are pairwise distinct. In particular, there exists an open neighborhood
$\mathcal{U}\subset B$ of $t_0$ on which the critical values
$\{u_i(t)\}_{i=1}^N$, $t\in \mathcal{U}$, form a {\em canonical} coordinate system, i.e., if we put
$1/\Delta_i:=(\partial/\partial u_i,\partial/\partial u_i)$ $(1\leq
i\leq N)$, then the map 
\ben
\Psi_t:\mathbb{C}^N\to T_t U,\quad e_i\mapsto
\sqrt{\Delta_i}\partial/\partial u_i
\een  
gives a trivialization of the tangent bundle $T\mathcal{U}$ in which the
Frobenius multiplication and the residue pairing have a diagonal form
\ben
e_i\bullet e_j = \delta_{i,j}\sqrt{\Delta_j} e_j,\quad (e_i,e_j)=\delta_{i,j}.
\een
If we write the differential equations \eqref{PF-1}--\eqref{de-prim} in
canonical coordinates, then we obtain a system of linear ODEs whose
coefficients depend polynomially on the entries of the matrices 
\ben
V:=-\Psi^{-1}\theta\Psi,
\quad 
(\lambda -U)^{-1},
\quad
\operatorname{ad}^{-1}_U(V), 
\een
where $U$ is the diagonal matrix
$\operatorname{diag}(u_1,\dots,u_N)$. Here
we are using that the matrix $V$ is skew-symmetric and we define
$\operatorname{ad}^{-1}_U(V)$ to be the unique matrix $X$ with
zero diagonal entries, s.t., $[U,X]=V$. The dependence of the matrix
$V$ on the canonical coordinates $\{u_i\}_{i=1}^N$ is quite
remarkable. It was proved by Dubrovin
that $V=V(u)$, $u=(u_1,\dots,u_N)$, is  a solution to an
integrable system (see equation (3.74) in \cite{Du}) and that $V$
extends analytically along any path in the domain $\mathfrak{D}=\{u\in 
\mathbb{C}^N\ |\ u_i\neq u_j,  \forall i\neq j\}$ (see Corollary 3.4 in \cite{Du}).   Using the canonical
coordinates, we embed $\mathcal{U}\subset \mathbb{C}^N$. Decreasing
$\mathcal{U}$ if necessary  we may assume that $\mathcal{U}\subset
\mathfrak{D}$. Let us pick a point $u^{(0)}\in \mathbb{C}^N$ whose coordinates are 
vertices of a regular polygon with center at $0$ and fix a path in
$\mathfrak{D}$ that connects $\mathcal{U}$
and $u^{(0)}$. We can extend $V(u)$ analytically along the path and obtain a
function analytic in a slightly larger domain
$\widetilde{\mathcal{U}}$. In particular, we can extend the system
\eqref{PF-1}--\eqref{de-prim} on $\widetilde{\mathcal{U}}\times
\mathbb{C}$, which gives us the required extension of the period
vectors. Again, we redenote by $B$ the extended domain obtained by
gluing $B$ and $\widetilde{\mathcal{U}}$ 
along the maximal open  subset $\mathcal{U}'$, s.t.,
$\mathcal{U}\subset \mathcal{U}'$ and the analytic embedding $\mathcal{U}\to
\widetilde{\mathcal{U}}$ extends to an analytic embedding $\mathcal{U}'\to
\widetilde{\mathcal{U}}$.

\subsection{The phase factors}
Our main goal is to investigate the analytic properties of the following infinite series:
\beq\label{inf-series}
\Omega_{\alpha,\beta}(t,\lambda,\mu) = \sum_{n=0}^\infty
(-1)^{n+1}(I^{(n)}_\alpha(t,\lambda),I^{(-n-1)}_\beta(t,\mu))\quad
\in\quad 
\mathbb{C}(\!(\lambda^{-1/|\sigma|})\!),
\eeq
where $\sigma$ is the classical monodromy operator and $|\sigma|$ is
the order of its semi-simple part. To begin with, it is not hard to
determine the radius of convergence. Namely, for fixed $(t,\mu)\in
S'$, the series \eqref{inf-series} is convergent for 
\ben
|\lambda|>\operatorname{max}_i\{|\mu|,|u_i(t)|\},
\een
where $u_i(t)$ is the set of critical values of $F(t,x)$. The main
motivation to study the above series comes from the so called {\em
  phase factors}
$B_{\alpha,\beta}(t,\lambda,\mu):=e^{\Omega_{\alpha,\beta}(t,\lambda,\mu)
}$, which are produced naturally when we compose two vertex operators
(see Section \ref{sec:VOA}).  
 
Recall the variation operator isomorphism
\ben
\operatorname{Var}:
\mathfrak{h}:=H^n(X_{0,\epsilon_0};\mathbb{Z})\to  
\mathfrak{h}^*:=H_n(X_{0,\epsilon_0};\mathbb{Z})
\een
defined by composing the Lefschetz duality
$H^n(X_{0,\epsilon_0};\mathbb{Z})\cong
H_n(X_{0,\epsilon_0},\partial X_{0,\epsilon_0}
;\mathbb{Z})$ and the operator $h_*-1$, where
$h:(X_{0,\epsilon_0},\partial X_{0,\epsilon_0})\to
(X_{0,\epsilon_0},\partial X_{0,\epsilon_0})$ is the geometric
monodromy. Let us recall also the {\em Seifert forms} (see \cite{AGV})
in respectively cohomology and homology:  
\ben
\operatorname{SF}(A,B)& := &  (-1)^{n/2+1} \langle A,
\operatorname{Var}(B)\rangle ,\quad A,B\in \mathfrak{h},\\
\operatorname{SF}(\alpha,\beta) & := & 
(-1)^{n/2+1} \langle \operatorname{Var}^{-1}( \alpha),
\beta\rangle,\quad \alpha,\beta\in \mathfrak{h}^*,
\een
where recall that we required $n$ to be even. Both forms are
non-degenerate, integer valued, and 
$\sigma$-invariant. The sign here is chosen in such a way that the
symmetrization of the Seifert form yields the sign-normalized
intersection form introduced above
\ben
\operatorname{SF}(\alpha,\beta) +\operatorname{SF}(\beta,\alpha) =(\alpha|\beta).
\een
Using the variation isomorphism we introduce intersection form in cohomology
$(A|B):=(\operatorname{Var}(A)|\operatorname{Var}(B))$. Note that the
cohomological intersection form is also the symmetrization of the
cohomological Seifert form. 

Finally, let us introduce the following polylogorithm operator
series depending on a linear operator $\sigma$ whose
eigenvalues are roots of unity:
\ben
\operatorname{Li}_\sigma(x) = 
\sum_{k=1}^\infty \frac{x^{k+\mathcal{N}}}{k+\mathcal{N}} \,,
\een
where 
$
\mathcal{N}:=-\frac{1}{2\pi\sqrt{-1}  } \log\sigma,
$
s.t., the eigenvalues of $\mathcal{N}$ belong to the set $(-1,0]\cap\mathbb{Q}$. 
Note that the series is convergent for $|x|<1$. Moreover, it can be expressed
in terms of the standard polylogorithm functions
\beq\label{poly}
\operatorname{Li}_\sigma(x) = x^{\mathcal{N}_n} \sum_{p=1}^\infty \sum_{r=1}^{|\sigma|}
(-\mathcal{N}_n|\sigma|)^{p-1}
\operatorname{Li}_p(\eta^r x^{1/|\sigma|}) \sigma_s^r,
\eeq
where we decomposed $\mathcal{N}=\mathcal{N}_s+\mathcal{N}_n$ into a
diagonal and nilpotent operator,
$\sigma_s:=e^{-2\pi\sqrt{-1}\mathcal{N}_s}$, $|\sigma|$ is the
order of $\sigma_s$, and $\eta=e^{2\pi\sqrt{-1}/|\sigma|}$. In particular, the series $\operatorname{Li}_\sigma(x)$
extends analytically along any path avoiding $x=0$ and $x=1$. 
\begin{theorem}\label{t1}
The following formula holds
\ben
\Omega_{\alpha,\beta}(0,\lambda,\mu) =  
-(\operatorname{Li}_\sigma(\mu/\lambda) \alpha|\beta)
+ P_{\alpha,\beta}(\lambda,\mu),
\een
where $P_{\alpha,\beta}(\lambda,\mu)$ is a polynomial in
$\lambda^{\pm 1/|\sigma|}, \mu^{\pm 1/|\sigma|} $,
$\log \lambda$, and $\log \mu$ satisfying 
\ben
P_{\alpha,\beta}(\lambda,\mu)-P_{\beta,\alpha}(\mu,\lambda) =
\operatorname{SF}\Big( 
\frac{e^{-2\pi\sqrt{-1}\mathcal{N}}-1}{\mathcal{N}}\, (\mu/\lambda)^\mathcal{N}\alpha_1,\beta\Big), 
\een
where $\alpha_1$ is the projection of $\alpha$ on the generalized
eigen subspace $\mathfrak{h}^*_1:=\operatorname{Ker}(\mathcal{N}_s).$
\end{theorem}
\noindent
As a corollary of Theorem \ref{t1} we obtain the following result (for
the proof see Lemma \ref{phase-sym:1}).
\begin{corollary}\label{c1}
Let $(\lambda,\mu)\in \mathbb{C}^2$ be a fixed point, s.t.,
$|\lambda|>|\mu|>0$ and $|\lambda-\mu|\ll 1$. Let $C\subset
\mathbb{C}^2$ be a path from $(\lambda,\mu)$ to $(\mu,\lambda)$
avoiding the diagonal and contained in a sufficiently small neighborhood of the
line segment $[(\lambda,\mu),(\mu,\lambda)]$.      
The following symmetry holds
$B_{\alpha,\beta}(0,\lambda,\mu)=B_{\beta,\alpha}(0,\mu,\lambda)$,
where the second phase factor is obtained from
$B_{\beta,\alpha}(0,\lambda,\mu)$ via analytic continuation along $C$.
\end{corollary}
Corollary \ref{c1} is important for the applications to VOA
representations, because it essentially means that the vertex
operators are mutually local. In particular, we can construct a
representation of the lattice VOA corresponding to the Milnor
lattice. In Section \ref{sec:VOA} we give more details on how to
construct such a representation, although a more systematic
investigation will be presented elsewhere.  

\subsection{The phase form and analytic extension}
The phase 1-form is by definition the following 1-form on $B'$
depending on a parameter $\xi$:
\beq\label{pf}
\mathcal{W}_{\alpha,\beta}(t,\xi) = I^{(0)}_\alpha(t,\xi)\bullet_t
I^{(0)}_\beta(t,0),\quad \alpha,\beta\in \mathfrak{h}^*,
\eeq
where the RHS is interpreted as a cotangent vector in $T_t^*B$ via the
identifications \eqref{flat-triv} and the parameter $\xi$ is assumed
to be sufficiently small. In other words, we expand the RHS into a
Taylor series at $\xi=0$
\ben
\sum_{m=0}^\infty
\mathcal{W}^{(m)}_{\alpha,\beta}(t)\frac{\xi^m}{m!},\quad
\mathcal{W}^{(m)}_{\alpha,\beta}(t) = I^{(m)}_\alpha(t,0)\bullet
I_\beta^{(0)}(t,0)\in T^*_tB.
\een
where for each $t\in B'$ the radius of convergence of the above series
is non-zero. 

The phase form can be used to describe the dependence of the phase
factors in the deformation parameters
\beq\label{phase-ac}
\Omega_{\alpha,\beta}(t,\lambda,\mu)=
\Omega_{\alpha,\beta}(0,\lambda,\mu) + 
\int_{-\lambda\mathbf{1}}^{t-\lambda\mathbf{1}} \mathcal{W}_{\alpha,\beta} (t',\mu-\lambda),
\eeq
where $(t,\lambda,\mu)$ is in the
domain of convergence of \eqref{inf-series}, $|\mu-\lambda|\ll 1$, and
the integration path should be chosen appropriately. 
The main result of this paper can be stated as follows. 
\begin{theorem}\label{t2}
Let $(t,\lambda,\mu)$ be a point in the domain of convergence of
\eqref{inf-series}, s.t., $|\mu-\lambda|\ll 1$ and $C\subset B'$ be a closed
loop based at $t-\lambda\mathbf{1}$, then 
\ben
\oint_{t'\in C} \mathcal{W}_{\alpha,\beta}(t',\mu-\lambda)
-\Omega_{w(\alpha),w(\beta)}(t,\lambda,\mu)
+\Omega_{\alpha,\beta}(t,\lambda,\mu)\quad \in \quad 2\pi\sqrt{-1}\mathbb{Z},
\een
where $w$ is the monodromy transformation along $C$. 
\end{theorem}
The integral \eqref{phase-ac} provides analytic extension along 
any path, so the phase factors can be interpreted as 
Laurent series in $(\mu-\lambda)$. Theorem \ref{t2} implies that
the coefficients of the Laurent series expansion are analytic on the
monodromy covering space, i.e, if a closed loop $C$ is in the kernel of
the monodromy representation, then the coefficients are invariant with
respect to the analytic continuation along $C$. Finally, when the
cycles $\alpha$ and $\beta$ are invariant along $C$, we get an
affirmative answer to Givental's question from \cite{G3}. 
\begin{corollary}\label{c2}
Let $C$ be a closed loop in $B'$ and $\alpha,\beta\in Q$ be cycles
invariant with respect to the monodromy transformation along $C$, then 
$\oint_{t\in C} \mathcal{W}_{\alpha,\beta} (t,\xi)\in 2\pi\sqrt{-1} \mathbb{Z}$.
\end{corollary}

\section{The phase factors at $t=0$}

Let $\mathcal{H}=H(\!(z^{-1})\!)$ be Givental's symplectic loop space,
where the symplectic form is 
\ben
\Omega(f(z),g(z))=\operatorname{Res}_{z=0} (f(-z),g(z))dz,\quad f,g\in \mathcal{H}.
\een
There is a natural polarization $\mathcal{H}=\mathcal{H}_+\oplus
\mathcal{H}_-$ with $\mathcal{H}_+:=H[z]$ and
$\mathcal{H}_-:=H[\![z^{-1}]\!]z^{-1}$ which allows us to identify
$\mathcal{H}\cong T^*\mathcal{H}_+$. 
Let us introduce the generating series  
\ben
\mathbf{f}_a(t,\lambda;z) = \sum_{n\in \mathbb{Z}}
I^{(n)}_a(t,\lambda)\, (-z)^n ,\quad a\in \mathfrak{h}^*
\een
and note that
$\Omega_{\alpha,\beta}(t,\lambda,\mu)=
\Omega(\mathbf{f}_\alpha(t,\lambda;z)_+,\mathbf{f}_\beta(t,\mu;z))$,
where the index $+$ denotes the projection on $\mathcal{H}_+$ along
$\mathcal{H}_-$.  

\subsection{The Virasoro grading operator}
The symplectic vector space $\mathcal{H}$ has the following Virasoro
grading operator:
\ben
\ell_0(z)=z\partial_z+\frac{1}{2}-\theta +\rho/z,
\een
where $\rho\in \operatorname{End}(H)$ is the operator of
multiplication by $W$. Note that $\rho/z$ is a nilpotent operator,
commuting with $\ell_0(z)$, while $\ell_0(z)-\rho/z$ is diagonalizable. 
Let us decompose, the symplectic vector space
\beq\label{energy}
\mathcal{H}=\mathcal{H}_{<0}\oplus \mathcal{H}_0\oplus \mathcal{H}_{>0},
\eeq
according to the sign of the eigenvalues of $\ell_0(z)$, i.e.,
$\mathcal{H}_{<0}$ is the sum of all generalized eigen subspaces of
$\ell_0$ with eigenvalue $<0$, $\mathcal{H}_{0}$ -- with eigenvalue 0,
and $\mathcal{H}_{>0}$ -- with eigenvalue $>0$. 

We are going to compute the symplectic pairing
\ben
\Omega(\mathbf{f}_\alpha(0,\lambda;z)_{>
  0},\mathbf{f}_\beta(0,\mu;z)) =
\Omega(\mathbf{f}_\alpha(0,\lambda;z),\mathbf{f}_\beta(0,\mu;z) _{<0} )
,\quad \alpha,\beta\in \mathfrak{h}^*, 
\een
where the index $> 0$ (resp. $<0$) corresponds to projection with respect to the
spectral decomposition of the operator $\ell_0(z)$ onto the subspace
spanned by all generalized eigen subspaces with eigenvalue $>0$
(resp. $<0$). It is not hard to check that the 
difference between the phase factor
$\Omega_{\alpha,\beta}(0,\lambda,\mu)$ and the above symplectic
pairing is a function $P_{\alpha,\beta}$ that has the form described
in Theorem \ref{t1}. So the main difficulty is in finding an explicit
formula for the above symplectic pairing.

\subsection{Fundamental solution}
The differential equation \eqref{de-prim} can be solved explicitly
when $t=0$. Assume that $k=-m$ with $m\gg 0$, then the fundamental
solution is given by the following operator-valued function
\ben
\Phi_m(\lambda) = e^{\rho\partial_\lambda\partial_m}\Big(
\frac{\lambda^{\theta+m-1/2}}{\Gamma(\theta+m+1/2)}\Big), 
\een
where the RHS is defined first for $m$ a complex number, so that
$\partial_m$ makes sense and the Gamma function is defined through its
Taylor's series expansion at $\theta=0$, so it makes sense to
substitute a linear operator for $\theta$. Note that if we assume that
$1/\Gamma(m)=0$ when $m$ is a negative integer or $0$, then the above formula
makes sense for all $m\in \mathbb{Z}$. 

Using the commutation relations
\beq\label{com-rel}
[\rho,\theta]=\rho,\quad \rho\partial_\lambda \Phi_m(\lambda) =
\Phi_m(\lambda)\rho, 
\eeq
it is easy to prove that the analytic continuation around $\lambda=0$
transforms $\Phi_m(\lambda)$ into 
\ben
\Phi_m(\lambda)\, e^{2\pi\sqrt{-1}\rho}\, e^{2\pi\sqrt{-1}(\theta-1/2)}.
\een
\subsection{Geometric sections}
If $\omega\in \Omega^{n+1}_{X/B}(X)$ is a
holomorphic form then $s(\omega):=\nabla_{\partial_\lambda} \int
d^{-1}\omega,$ where $\nabla$ is the Gauss--Manin connection, defines
a holomorphic section of the vanishing 
cohomology bundle known as {\em geometric section}. Let us denote by
$s(\omega,\lambda)$ the value of the section at the point
$(0,\lambda)\in S'$. 

Put $\ell=n/2$ (we assume that $n$ is even). Note that by definition the period vector 
\ben
I^{(-\ell)}(0,\lambda) = (2\pi)^{-\ell} \sum_{i=1}^N \phi^i\otimes s(\omega_i,\lambda),
\een 
where $\{\phi^i\}\subset H$ is a basis of $H$ dual to $\{\phi_i\}$
with respect to the residue pairing and $\omega_i\in
\Omega^{n+1}_{X/B}(X)$ is defined by 
\ben
-\nabla_{\partial_i} \int d^{-1}\omega = \nabla_{\partial_\lambda}
\int d^{-1}\omega_i. 
\een
On the other hand, since $\Phi_\ell(\lambda)$ is a fundamental
solution for the differential equation \eqref{de-prim} with $k=-\ell$,
we get that 
\ben
I^{(-\ell)}(0,\lambda) = \sum_{i=1}^N \Phi_\ell(\lambda)\phi^i\otimes A_i,
\een
where $\{A_i\}_{i=1}^N\subset \mathfrak{h}$ is a basis with some very
special properties that guarantee that $\omega$ is a primitive form
(see \cite{He} for more details). 
 Let us introduce the linear map
\ben
{}^\#:\operatorname{End}(H)\to \operatorname{End}(\mathfrak{h}),
\quad 
\mbox{s.t., }\quad 
\sum_{i=1}^N T\phi^i\otimes A_i=\sum_{i=1}^N \phi^i\otimes T^\#A_i.
\een
\begin{lemma}\label{anti-hom}
The map ${}^\#$ is anti-homomorphism of algebras, i.e., $(T_1T_2)^\# = T_2^\#T_1^\#.$
\end{lemma}   
\proof
We just have to use that $T^\#A_i=\sum_{j=1}^N (T\phi^j,\phi_i)A_j$. \qed

The numbers 
\ben
s_i:=d_i+\ell-\frac{1}{2}-\frac{D}{2},\quad 1\leq i\leq N,
\een
are known as the {\em Steenbrink spectrum} of the singularity. The key
result that reveals the Hodge-theoretic origin of the period integrals
is the following Lemma.
\begin{proposition}\label{geom-sec}
The following formula holds
\ben
s(\omega_i,\lambda) = (2\pi)^\ell\,\lambda^{s_i+\rho^\#} \Gamma(s_i+\rho^\#+1)^{-1}A_i,
\een
where the value of the Gamma function is defined through its Taylor's
series at $\rho^\#=0$. 
\end{proposition}  
\proof
We have to compute $\Phi_\ell(\lambda)^\#A_i$. Using the commutation
relations \eqref{com-rel} we get
\ben
\Phi_\ell(\lambda) = \sum_{k=0}^\infty
\frac{1}{k!}\rho^k\partial_\lambda^k\partial_\ell^k \Big(
\lambda^{\theta+\ell-1/2}\Gamma(\theta+\ell+1/2)^{-1}\Big) = 
\sum_{k=0}^\infty
\frac{1}{k!}\partial_\ell^k \Big(
\lambda^{\theta+\ell-1/2}\Gamma(\theta+\ell+1/2)^{-1}\Big) \rho^k.
\een
Note that 
\ben
\theta^\# A_i = \sum_{j=1}^N (\theta\phi^j,\phi_i)A_j = (d_i-D/2)A_i. 
\een
To finish the proof we just need to recall Lemma \ref{anti-hom} and
use Taylor's formula.
\qed

Since the geometric sections are single valued, Proposition
\ref{geom-sec} allows us to conclude that the classical monodromy 
\ben
\sigma = e^{-2\pi\sqrt{-1} \rho^\#}\, e^{-2\pi\sqrt{-1}(\theta^\#-1/2)}.
\een
Let us write $s_i=p_i+\alpha_i$, where $-1<\alpha_i\leq 0$, then our
choice of logarithm  $\mathcal{N}=-\frac{1}{2\pi\sqrt{-1}} \log
\sigma$ acts on the basis $\{A_i\}\subset \mathfrak{h}$ as follows
\ben
\mathcal{N}(A_i) = (\alpha_i +\rho^\#) A_i,\quad 1\leq i\leq N.
\een
Let us point out that since
\ben
\nabla^{p_i}_{\partial_\lambda}s(\omega_i,\lambda) =
\lambda^{\alpha_i+\rho^\#} \Gamma(\alpha_i+\rho^\#+1)^{-1} A_i
\een
the vectors $\Gamma(\alpha_i+\rho^\#+1)^{-1} A_i\in
F_{p_i}\mathfrak{h}$, where $\{F_p\mathfrak{h}\}_{p=0}^n$ is the Steenbrink's
Hodge filtration.

\subsection{The higher residue pairing}

Let us denote by $\mathbb{H}\to B$ the vector bundle whose fiber over
$t\in B$ is
\ben
\mathbb{H}_t = \Omega_{X_t}^{n+1}[\![z]\!]/(z d_x + d_xF(t,x)\wedge) \Omega_{X_t}^{n}[\![z]\!],
\een
where
$X_t=X\cap \mathbb{C}^{n+1}\times \{t\}$. By definition, the domain
$X$ is chosen so small that 
$\mathbb{H}_t$ is a $\mathbb{C}[\![z]\!]$-free module of rank
$N$. We can think of sections $\omega$ of $\mathbb{H}$ as formal oscillatory
integrals $\int e^{F/z}\omega$, which allows us to see that
$\mathbb{H}$ is equipped with a flat connection, called {\em
  Gauss--Manin connection}, corresponding to
differentiating the integral formally with respect to the deformation
parameters. 

Recall, also K. Saito's higher residue pairing
\ben
K_t:\mathbb{H}_t\otimes \mathbb{H}_t\to \mathbb{C}[\![z]\!]z^{n+1}.
\een  
It is uniquely determined up to a constant by the following properties
\begin{enumerate}
\item[(K1)] If $\omega_1,\omega_2\in \mathbb{H}_t$, then
$K_t(\omega_1,\omega_2) = (-1)^{n+1}K_t(\omega_2,\omega_1)^*,$
where $*$ is the involution $z\mapsto -z$.
\item[(K2)]
If $p(z)\in \mathbb{C}[z]$, then
\ben
p(z)K_t(\omega_1,\omega_2) =K_t(p(z)\omega_1,\omega_2) =K_t(\omega_1,p(-z)\omega_2). 
\een
\item[(K3)]
The pairing $K^{(p)}_t(\omega_1,\omega_2)$ defined by the coefficient in front of
$z^{n+1+p}$ in $K_t(\omega_1,\omega_2)$ depends analytically on $t$
and the Leibnitz rule holds
\ben
\xi K_t(\omega_1,\omega_2) = K_t(\nabla_\xi\omega_1,\omega_2)+K_t(\omega_1,\nabla_\xi\omega_2),
\een
where $\xi=\partial/\partial t_i$ or $z\partial_z$ and  $\nabla$ is the
Gauss--Manin connection on $\mathbb{H}$.
\item[(K4)]
If $\omega_i=\phi_i(x)dx_0\cdots dx_n$ $(i=1,2)$, then
$K^{(0)}_t(\omega_1,\omega_2) = (\phi_1,\phi_2)$, where the residue
pairing (see \eqref{res-pairing}) is with respect to the volume form
$dx_0\cdots dx_n. $
\end{enumerate}

Given a holomorphic form $\omega\in \Omega^{n+1}_{X_0}(X_0)$, let us
denote by 
\ben
\widehat{s}(\omega,z) = (-2\pi z)^{-\ell-1/2} \int_0^\infty e^{\lambda/z} s(\omega,\lambda)d\lambda,
\een
where the integration path is $\lambda=-t z,$ $t\in [0,+\infty)$, the
Laplace transform of the corresponding geometric section. We will make
use also of the automorphism
$\widehat{s}(\omega,z)^*:=\widehat{s}(\omega,e^{\pi\sqrt{-1}}z).$ Note
that since we might have $\log z$ dependence, ${}^*$ is no longer an
involution.

Let us introduce also the so called {\em elementary sections}. They are
defined as follows. Given a vector $A\in \mathfrak{h}$, we can construct a section of the
vanishing cohomology bundle over $S'\cap\{t=0\}$ as follows
\ben
s(A,\lambda):=\sum_{-1<\alpha\leq 0} \lambda^{\alpha+\rho^\#} A_\alpha, 
\een
where $A_\alpha$ is the projection of $A$ on the generalized
eigen space of $\sigma$ corresponding to the eigenvalue
$e^{-2\pi\sqrt{-1}\alpha}$. Put
\ben
\widehat{s}(A,z) = (-2\pi z)^{-\ell-1/2} \int_0^\infty e^{\lambda/z} s(A,\lambda)d\lambda,
\een
where the integration path is the same as above. The next result is a
reformulation of Hertling's formula for Saito's higher residue
pairing in terms of the Seifert form.
\begin{lemma}\label{He-formula}
The pairing
\ben
K_W(\omega_1,\omega_2):=
\operatorname{SF}(\widehat{s}(\omega_1,z)^*,\widehat{s}(\omega_2,z))^*\, z^{n+1}
\een
coincides with Saito's higher residue pairing
$K_0(\omega_1,\omega_2)$. 
\end{lemma}
\proof
We will prove that the Lemma is equivalent to Hertling's formula for
Saito's higher residue pairings. In order to compare our formula to
Hertling's one, we have to introduce the vector bundle
$\mathbb{H}''\to B$ whose fiber over a point $t\in B$ is 
$\mathbb{H}''_t:=\Omega^{n+1}_{X_t}/d_xF(t,x)\wedge d_x
\Omega^{n-1}_{X_t}$. If $\omega\in \Omega^{n+1}_{X_t} $ is a
holomorphic form, then we denote by $s(\omega)\in \mathbb{H}''_t$ the projection of $\omega$ . It is known that $\mathbb{H}''_t$ is a
$\mathbb{C}[\partial_\lambda^{-1}]$-module and that if we identify
$\mathbb{C}[\![z]\!]\cong \mathbb{C}[\![\partial_\lambda^{-1}]\!]$ via
$z\mapsto -\partial_{\lambda}^{-1}$, then the natural map 
\ben
\mathbb{H}_t\to \mathbb{H}_t''\otimes_{\mathbb{C}[\partial_\lambda^{-1}]}
\mathbb{C}[\![\partial_\lambda^{-1}]\!] ,\quad 
\omega\mapsto s(\omega),
\een 
is an isomorphism of $\mathbb{C}[\![z]\!]$-modules. Let us define
\ben
\widetilde{K}_t(s(\omega_1),s(\omega_2)) :=
K_t(\omega_2,\omega_1)=(-1)^{n+1}K_t(\omega_1,\omega_2)^*.  
\een
It is easy to check that $\widetilde{K}$ satisfy properties (K1)--(K4)
of the higher residue pairing except that we have to replace $z$ by
$\partial_\lambda^{-1}$. Since the higher residue pairing is uniquely
determined by its properties, $\widetilde{K}$ coincides with the
higher residue pairing used in \cite{He}. We have to prove that 
\ben
\widetilde{K}_t(s(\omega_1),s(\omega_2)) =
\operatorname{SF}(\widehat{s}(\omega_1,z)^*,
\widehat{s}(\omega_2,z)) z^{n+1}.
\een
Note that using the above formula we can uniquely define a pairing on the space
of all elementary sections, s.t., property (K2) holds. We just need to check that the RHS agrees
with Hertling's formula for a pair of elementary sections (see
formulas (10.81)--(10.83) in \cite{He}). 

Let us assume that $B_i$ $(i=1,2)$ are eigenvectors with eigenvalues 
$e^{-2\pi\sqrt{-1}\beta_i}$, $-1<\beta_i\leq 0$. Then 
\ben
\widehat{s}(B_i,z) = (2\pi)^{-\ell-1/2}
(-z)^{\beta_i+\rho^\#-\ell+1/2} \Gamma(\beta_i+\rho^\#+1) B_i,\quad i=1,2.
\een
Using that the Seifert form is monodromy invariant and infinitesimally 
$\rho^\#$-invariant, i.e.,
$\operatorname{SF}(\rho^\#A,B)+\operatorname{SF}(A,\rho^\#B)=0$ 
we get that
\beq\label{SF-pairing}
\operatorname{SF}(\widehat{s}(B_1,z)^*,\widehat{s}(B_2,z)) \, z^{n+1}
\eeq
is given by the following formula
\ben
\frac{z^{\beta_1+\beta_2+2}}{(2\pi\sqrt{-1})^{n+1}}
\operatorname{SF}(B_1,
e^{\pi\sqrt{-1}(\beta_2+\rho^\#+1)}\Gamma(\beta_1-\rho^\#+1)\Gamma(\beta_2+\rho^\#+1)B_2). 
\een
Since the Seifert form is monodromy invariant, the above pairing
vanishes unless $\beta_1+\beta_2\in \mathbb{Z}$, i.e.,
$\beta_1+\beta_2=-1$, or $\beta_1=\beta_2=0.$ Recall, the operator
$\mathcal{N}=-\frac{1}{2\pi\sqrt{-1}} \log \sigma$ and note that
$\mathcal{N} B_i= - (\beta_i +\rho^\#)B_i$.  
\ben
e^{\pi\sqrt{-1}(\beta_2+\rho^\#+1)}\Gamma(\beta_1-\rho^\#+1)\Gamma(\beta_2+\rho^\#+1)B_2
= 
\begin{cases}
-2\pi\sqrt{-1} (\sigma-1)^{-1} B_2 & \mbox{ if } \beta_1+\beta_2=-1,\\
2\pi\sqrt{-1}  \mathcal{N} (\sigma-1)^{-1} B_2 & \mbox{ if } \beta_1=\beta_2=0.
\end{cases}
\een
Hence the pairing \eqref{SF-pairing} takes the form
\ben
\frac{1}{(2\pi\sqrt{-1})^n}\, (-1)^\ell  \langle B_1,\operatorname{Var}
\, (\sigma-1)^{-1} B_2\rangle \, z, \quad \mbox{ for } \beta_1+\beta_2=-1,
\een
\ben
\frac{-1}{(2\pi\sqrt{-1})^{n+1}} \, (-1)^\ell \langle B_1,\operatorname{Var}
\, 2\pi\sqrt{-1}  \mathcal{N} (\sigma-1)^{-1} B_2\rangle\, z^2, \quad \mbox{ for }
\beta_1=\beta_2=0,
\een
and it vanishes in all other cases. \qed

The identity in Lemma \ref{He-formula} is really remarkable. 
It is a relation between two completely different quantities. The LHS
is defined via residues of differential forms, while the RHS is purely
topological. Let us outline a different way to prove Lemma
\ref{He-formula}, which in particular generalizes the identity for an
arbitrary deformation. Given $\omega\in \Omega^{n+1}_{X/B}(X)$ and a
critical value $u=u(t)$ of $F(t,x)$ we define the following
formal asymptotic series
\ben
\widehat{s}_u(\omega,z) := (-2\pi z)^{-\ell-1/2} \int_{u(t)}^\infty
e^{\lambda/z} s(\omega) d\lambda.
\een 
Note that this is the stationary phase asymptotic as $z\to 0$ of an appropriate
oscillatory integral. Let us define
\ben
K_F(\omega_1,\omega_2) = \sum_{u}
\operatorname{SF}(\widehat{s}_u(\omega_1,z)^*,
\widehat{s}_u(\omega_1,z))^*\, z^{n+1},
\een
where the sum is over all critical values of $F$. When $F=W$, this
formula reduces to Lemma \ref{He-formula}. One has to check
that the above formula satisfies all properties (K1)--(K4) of the
higher residue pairing. The verification is straightforward except for
property (K4). In the latter case we take a generic deformation and
then we just have to see that the computation on the RHS reduces
to proving Lemma \ref{He-formula} for $A_1$-singularity, which is
straightforward.   It remains only to recall the result of M. Saito
\cite{MSa} that the higher residue 
pairing is uniquely determined by its properties. 

\subsection{Proof of Theorem \ref{t1}}
Using the explicit formula for the geometric sections in Proposition
\ref{geom-sec} we get that 
\ben
\widehat{s}(\omega_i,z) = (2\pi)^{-1/2} (-z)^{p_i+\mathcal{N}-\ell+1/2}A_i,
\een
where $p_i=\lfloor s_i\rfloor $ is the floor of the Steenbrink number
$s_i$. It is convenient, to introduce the operator $p\in
\operatorname{End}(\mathfrak{h})$, s.t., 
$p(A_i)=p_iA_i$.
Recalling Lemma \ref{He-formula} we get
\beq\label{res-sf}
K_W(\omega_i,\omega_j)=
(\phi_i,\phi_j) z^{n+1}= \frac{1}{2\pi\sqrt{-1}} \langle A_i, e^{\pi\sqrt{-1}
  \mathcal{N}} e^{\pi\sqrt{-1} p} A_j\rangle\, z^{s_i+s_j-2\ell+1 +n+1},
\eeq
where $\langle \, , \, \rangle$ is the Seifert form (without the sign
normalization) and the 1st identity holds, because $\omega$ is a
primitive form. Let us define the {\em residue pairing} $(\, ,\, )$ on
$\mathfrak{h}$, s.t., the $(\phi_i,\phi_j)=(A_i,A_j)$, i.e., 
\ben
(A,B) = \frac{1}{2\pi\sqrt{-1}} \langle A, e^{\pi\sqrt{-1}
  \mathcal{N}} e^{\pi\sqrt{-1} p}B\rangle,\quad 
A,B\in \mathfrak{h}.
\een 
Finally, if
$R\in \operatorname{End}(\mathfrak{h})$, then we denote by $R^T$ the
transpose with respect to the residue 
pairing. Let us point out the following relations
\ben
\mathcal{N}=\mathcal{N}_s+\rho^\#,\quad p+\mathcal{N}_s= (\theta +\ell-1/2)^\#,\quad [\mathcal{N},\rho^\#]=0,
\quad
[p,\rho^\#]=\rho^\#,
\een
where $\mathcal{N}_s$ is the semi-simple part of
$\mathcal{N}$. Comparing the powers of $z$  in \eqref{res-sf} we get
\ben
p+\mathcal{N}_s+p^{T}+\mathcal{N}_s^{T} = 2\ell-1.
\een
By definition 
\ben
\mathbf{f}_\alpha(0,\lambda;z) = \sum_{m\in
  \mathbb{Z}}\partial_\lambda^{m+\ell}
I^{(-\ell)}_\alpha(0,\lambda)\,(-z)^m =
\sum_{i=1}^N \sum_{m\in  \mathbb{Z}} 
(-z)^m \partial_\lambda^{m+\ell-p_i}
\langle\lambda^{\mathcal{N}}\Gamma(\mathcal{N}+1)^{-1} A_i,\alpha\rangle \phi^i
\een
Let us introduce also the transpose  $R^{SF}$ with respect
to the Seifert form
\ben
\langle RA,B\rangle =: \langle A, R^{SF}B\rangle,\quad \forall A,B\in \mathfrak{h}.
\een 
Note that since the Seifert form is not symmetric, in general 
$ \langle A,RB\rangle\neq \langle R^{SF}A,B\rangle$.
Since 
\ben
\Big(z\partial_z+\frac{1}{2}-\theta\Big) (-z)^m\phi^i =  (m+\ell-s_i)  (-z)^m\phi^i,
\een
after shifting $m\mapsto m-\ell+p_i$ we get that the projection 
\beq\label{proj<0}
\mathbf{f}_\beta(0,\mu;z)_{<0}=
\sum_{j=1}^N \sum_{m=-\infty}^{-1} 
\left\langle
A_j,
(-z)^{m-\ell+p^{SF}}
\mu^{\mathcal{N}^{SF}-m}\Gamma(\mathcal{N}^{SF}-m+1)^{-1}
\beta
\right\rangle\phi^j .
\eeq
Similarly,
\ben
\mathbf{f}_\alpha(0,\lambda;z)=
\sum_{i=1}^N \sum_{k\in \mathbb{Z}}
\left\langle
A_i, 
(-z)^{k-\ell+p^{SF}}
\lambda^{\mathcal{N}^{SF}-k}\Gamma(\mathcal{N}^{SF}-k+1)^{-1}
\alpha\right\rangle \phi^i .
\een
Let us denote by $\{A^i\}_{i=1}^N$ the basis of $\mathfrak{h}$ dual to
$\{A_i\}_{i=1}^N$ with respect to the residue pairing. We need to
compute a pairing of the following type:
\ben
\sum_{i,j=1}^N \langle A_i,x\rangle (\phi^i,\phi^j) \langle 
A_j,y\rangle =
(2\pi\sqrt{-1})^2\sum_{i,j=1}^N (A_i, M^{-1} x) (\phi^i,\phi^j)
(A_j, M^{-1} y) ,
\een
where $M=e^{\pi\sqrt{-1}\mathcal{N}}e^{\pi\sqrt{-1}p}$ and the
equality follows from the definition of the residue pairing.
Using that $(\phi^i,\phi^j)=(A^i,A^j)$ and the standard properties of
dual bases we get 
\ben
(2\pi\sqrt{-1})^2 ((M^{-1}) x, (M^{-1}) y) = 
2\pi\sqrt{-1}\langle M^{-1}x,y\rangle.
\een
The symplectic pairing 
\begin{align}\label{vir-prop}
&&
\Omega(\mathbf{f}_\alpha(0,\lambda;z)_{>0}, \mathbf{f}_\beta(0,\mu;z))=
\Omega(\mathbf{f}_\alpha(0,\lambda;z),\mathbf{f}_\beta(0,\mu;z)_{<0})
 = 
2\pi\sqrt{-1}\,
\sum_{m=-\infty}^{-1} \sum_{k\in \mathbb{Z}} 
\operatorname{Res}_{z=0} \\
\notag
&&
\left\langle M^{-1}
z^{k-\ell+p^{SF}}
\lambda^{\mathcal{N}^{SF}-k}\Gamma(\mathcal{N}^{SF}-k+1)^{-1}
\alpha,
(-z)^{m-\ell+p^{SF}}
\mu^{\mathcal{N}^{SF}-m}\Gamma(\mathcal{N}^{SF}-m+1)^{-1}
\beta
 \right\rangle. 
\end{align}
Using that by definition $R^{SF}=MR^TM^{-1}$ we get the following
formulas
\ben
\mathcal{N}_s^T=\mathcal{N}_s^{SF},\quad
(\rho^\#)^{SF} = -\rho^\#,\quad
(\rho^\#)^T= \rho^\#,\quad 
p^{SF} = p^T+\pi\sqrt{-1}\rho^\#. 
\een
Let us point out that the two conjugations do not commute in
general and also that $(R^{SF})^{SF}\neq R$ in general. 
Since $M^{-1} p^{SF} = p^T M^{-1}$ and $M^{-1}
\mathcal{N}^{SF} = \mathcal{N}^T M^{-1}$, the summand in
the sum \eqref{vir-prop} can be written in the following form
\ben
\left\langle
z^{m+k-2\ell+p+p^T} e^{\pi\sqrt{-1} (m-\ell+p)}
\lambda^{\mathcal{N}^T-k}\Gamma(\mathcal{N}^T-k+1)^{-1}
M^{-1}\alpha,
\mu^{\mathcal{N}^{SF}-m}\Gamma(\mathcal{N}^{SF}-m+1)^{-1}
\beta
 \right\rangle. 
\een
Using that $e^{\pi\sqrt{-1} p}\mathcal{N}^T=\mathcal{N}^{SF}
e^{\pi\sqrt{-1} p}$ and 
$
e^{\pi\sqrt{-1} p} M^{-1} = e^{-\pi\sqrt{-1}\mathcal{N}},
$
we get
\ben
\left\langle
z^{m+k-2\ell+p+p^T} 
\lambda^{\mathcal{N}^{SF}-k}\Gamma(\mathcal{N}^{SF}-k+1)^{-1}
e^{-\pi\sqrt{-1}(\mathcal{N}-m+\ell)}
\alpha,
\mu^{\mathcal{N}^{SF}-m}\Gamma(\mathcal{N}^{SF}-m+1)^{-1}
\beta
 \right\rangle. 
\een 
The terms that contribute to the residue must satisfy 
\ben
m+k+p+p^T=2\ell-1 = p+p^T+\mathcal{N}_s+\mathcal{N}_s^T=
p+p^T+\mathcal{N}+\mathcal{N}^{SF},
\een
so we may substitute $p+p^T= 2\ell-1-m-k$ and $\mathcal{N}^{SF}-k =
-\mathcal{N}+m$
\ben
(-1)^{\ell}z^{-1}\left\langle
(\mu/\lambda)^{\mathcal{N}-m}
\Gamma(\mathcal{N}-m+1)^{-1}
\Gamma(-\mathcal{N}+m+1)^{-1}
e^{\pi\sqrt{-1}(-\mathcal{N}+m)}
\alpha, \beta
 \right\rangle. 
\een 
Using the well known identity
\ben
\frac{e^{\pi\sqrt{-1} x}}{\Gamma(1-x)\Gamma(1+x)} =
  \frac{e^{2\pi\sqrt{-1} x}-1}{2\pi\sqrt{-1} x}
\een
with $x=-\mathcal{N}+m$ we get
\ben
\frac{(-1)^{\ell+1}}{2\pi\sqrt{-1}}\,z^{-1}\,
\left\langle 
\frac{(\mu/\lambda)^{\mathcal{N}-m}}{\mathcal{N}-m}
(e^{-2\pi\sqrt{-1}\mathcal{N}}-1)
\alpha, \beta
\right\rangle.
\een
It remains only to recall that the intersection form can be expressed
in terms of the Seifert form  
\ben
(\alpha|\beta) = (-1)^\ell \langle (\sigma-1)\alpha,\beta
\rangle,\quad \alpha,\beta\in \mathfrak{h}^*. 
\een
We get
\ben
\Omega(\mathbf{f}_\alpha(0,\lambda;z)_{>0}, \mathbf{f}_\beta(0,\mu;z))
= -(\operatorname{Li}_\sigma(\mu/\lambda)\alpha|\beta).
\een
Put 
\ben
P_{\alpha,\beta}(\lambda,\mu) = 
\Omega(\mathbf{f}_\alpha(0,\lambda;z)_+, \mathbf{f}_\beta(0,\mu;z))-
\Omega(\mathbf{f}_\alpha(0,\lambda;z)_{>0}, \mathbf{f}_\beta(0,\mu;z)).
\een
Using that $\Omega$ is a symplectic pairing and that
$\mathcal{H}_{>0}$ is symplectic orthogonal to $\mathcal{H}_{\geq 0}$
we get 
\ben
P_{\beta,\alpha}(\mu,\lambda)=
\Omega(\mathbf{f}_\beta(0,\mu;z),\mathbf{f}_\alpha(0,\lambda;z)_-) - 
\Omega(\mathbf{f}_\beta(0,\mu;z), \mathbf{f}_\alpha(0,\lambda;z)_{<0}).
\een
Hence
\ben
P_{\alpha,\beta}(\lambda,\mu)-P_{\beta,\alpha}(\mu,\lambda) = 
\Omega(\mathbf{f}_\alpha(0,\lambda;z)_0,\mathbf{f}_\beta(0,\mu;z)_0).
\een
The above symplectic pairing  can be computed in the same way as
above. The only difference is that the summation over $m$ in
\eqref{proj<0} collapses to $m=0$ and the
summation over $j$ reduces only to $j$, s.t., $s_j=p_j$, i.e., $\mathcal{N}_s
A_j=0$. This implies that in \eqref{vir-prop} we have to put $m=0$ and
replace $\alpha$ (and $\beta$) by $\alpha_1$ (and $\beta_1$). The
rest of the computation is the same. We get
\ben
\Omega(\mathbf{f}_\alpha(0,\lambda;z)_0,\mathbf{f}_\beta(0,\mu;z)_0) = 
(-1)^{\ell+1}\,
\left\langle 
\frac{(\mu/\lambda)^{\mathcal{N}}}{\mathcal{N}}
(e^{-2\pi\sqrt{-1}\mathcal{N}}-1)
\alpha_1, \beta_1
\right\rangle.
\een
Finally, the
polynomiality statement about $P_{\alpha,\beta}(\lambda,\mu)$ follows
from the fact that the vector space $\mathcal{H}_+\cap
\mathcal{H}_{<0}$ is finite dimensional.\qed

\section{Analytic continuation} 

We will be interested in the phase factors as multivalued analytic
functions. Our starting point is the infinite series
\eqref{inf-series}, which by definition is interpreted as a formal
Laurent series of the type
\beq\label{laurent-series}
\sum_{r=0}^{|\sigma|-1} \sum_{\ell=0}^{d-1} \lambda^{m+r/|\sigma|} (\log \lambda)^{\ell}
\sum_{k=0}^\infty \Omega_{\alpha,\beta}^{r,\ell,k}(t,\mu) \lambda^{-k},
\eeq
where $d\geq 1$ is the smallest integer number, s.t., $\mathcal{N}_n^d=0$. If $R>0$
is the smallest real number, s.t., for every $r$ and $\ell$ the infinite series in
$k$ is convergent for all $|\lambda|>R$, then $R^{-1}$ is called the {\em radius of
  convergence}. The radius of convergence does not change if we
differentiate the series with respect to $\lambda$. On the other hand,
the derivative of \eqref{inf-series} with respect to $\lambda$ can be computed in terms of the
period integrals (see \cite{M3} Proposition 2.3)
\beq\label{pf-partial}
\partial_\lambda \Omega(\mathbf{f}_\alpha(t,\lambda;z)_+,
\mathbf{f}_\beta(t,\mu;z)) = \frac{1}{\lambda-\mu} 
\Big( I^{(0)}_\alpha(t,\lambda),(\theta+1/2)  I^{(-1)}_\beta(t,\mu)\Big).
\eeq
Since the period $I^{(0)}_\alpha(t,\lambda)$ is a solution to a
Fuchsian differential equation in $\lambda$, whose singularities are
at $\lambda=\infty$ and the critical values $u_i(t)$ of $F(x,t)$, we
get that the series \eqref{laurent-series} is convergent for 
\ben
|\lambda|>\operatorname{max}\{ |\mu|, r(t)\},
\een
where $r(t)=\operatorname{max}_i\{|u_i(t)|\}$.  
\subsection{The domains at infinity}
Put $S_\infty=\{(t,\lambda)\in S \, |\, |\lambda|>r(t)\}$. Let us fix
a real constant $c\in (0,1)$, say
$c=|1-e^{2\pi\sqrt{-1}/|\sigma|}|/3$,  s.t., the distance between 
any two distinct $|\sigma|$-roots of 1 is bigger than $c$. Let 
$
\epsilon: S_\infty \to \mathbb{R}_{>0}
$
be the function 
\ben
\epsilon(t,\lambda) =
\operatorname{sup}_{C}\Big(\operatorname{min}_{(t',\lambda')\in C}
(|\lambda'|-r(t')) \Big),
\een
where the $\operatorname{sup}$ is over all paths $C\subset S_\infty$
from $(0,\lambda)$ to $(t,\lambda)$.  
We introduce the domain  
\ben
D_\infty:=\Big\{(t,\lambda,\mu) \in S_\infty\times \mathbb{C}\ |\ 
|\lambda-\mu|<
c\,\operatorname{min}\{\epsilon(t,\lambda),\epsilon(t,\mu)\} \Big\} 
\een
and its subdomain 
$
D^+_\infty:=\{ |\lambda|>|\mu|\}\subset D_\infty.
$

Let us denote by $\widetilde{S}_\infty$ the
universal covering of $S_\infty$. The domains $D^+_\infty$ and
$D_\infty$ are (trivial) smooth disk fibrations over $S_\infty$. The
phase factors $\Omega_{\alpha,\beta}$ are holomorphic functions on the
pullback $\widetilde{D}_\infty^+$ of $D_\infty^+$ via the covering map
$\widetilde{S}_\infty\to S_\infty$. In more explicit terms, in order
to define the phase factor $\Omega_{\alpha,\beta}(t,\lambda,\mu)$ at
some point $(t,\lambda,\mu)\in D_\infty^+$ we have to select a
reference path in $S_\infty$ from $(0,\epsilon_0)$ to
$(t,\lambda)$. This choice determines the value of the period vectors
$I^{(n)}_\alpha(t,\lambda)$ and using the line segment
$[(t,\lambda),(t,\mu)]$ we can specify also the values of
$I^{(-n-1)}_\beta(t,\mu)$, so the summands in \eqref{inf-series} are
uniquely defined. Note that if $(t,\lambda,\mu)\in D_\infty$, then
$|\lambda-\mu|<|\lambda|-r(t)$, so using the triangle inequality it is
easy to verify that the line segment $[(t,\lambda),(t,\mu)]$  does
not intersect the discriminant. 

\subsection{Symmetry of the phase factors at $t=0$}
According to Theorem \ref{t1}, the phase factor
$\Omega_{\alpha,\beta}(0,\lambda,\mu)$ extends analytically from
$D_\infty^+$ to $D_\infty-\{\lambda=\mu\}$. Note that the
polylogorithms that enter in our formula have singularities for
$(\lambda,\mu)$, s.t., $\mu=0$ or $\lambda^{|\sigma|}=\mu^{|\sigma|}$.
However, if $(0,\lambda,\mu)\in D_\infty$, then $\mu\neq 0$ and thanks
to our choice of the constant $c$ the equality 
$\lambda^{|\sigma|} =\mu^{|\sigma|}$ implies that $\lambda =\mu$ . Therefore,
our explicit formula provides the analytic extension along any path
in $D_\infty -\{\lambda=\mu\}$.  The analytic continuation has the
following crucial symmetry. Fix $(0,\lambda,\mu)\in D_\infty^+$ and a
path $C\subset D_\infty -\{\lambda=\mu\}$ connecting $(0,\lambda,\mu)$
with $(0,\mu,\lambda)$. The coordinate projections of $C$ along the
last two coordinates determine paths $C_1$ from $(0,\lambda)$ to
$(0,\mu)$ and $C_2$ from $(0,\mu)$ to $(0,\lambda)$ in $S_\infty$.  
\begin{lemma}\label{phase-sym:1}
Let $C\subset D_\infty-\{\lambda=\mu\}$ be a path from
$(0,\lambda,\mu)\in D_\infty^+$ to $(0,\mu,\lambda)$, s.t., the projections $C_1$
and $C_2$ of $C$ are homotopic (in $S_\infty$) respectively to the line segments
$[(0,\lambda),(0,\mu)]$ and $[(0,\mu),(0,\lambda)]$,  then 
\ben
\Omega_{\alpha,\beta}(0,\lambda,\mu)-\Omega_{\beta,\alpha}(0,\mu,\lambda)
=- 2\pi\sqrt{-1} \Big(
\operatorname{SF}(\alpha,\beta)+k(\alpha|\beta)
\Big),\quad k\in \mathbb{Z},
\een 
where the 2nd phase factor is obtained from
$\Omega_{\beta,\alpha}(0,\lambda,\mu)$ via analytic continuation along
the path $C$ and the integer $k$ depends on the choice of $C$.
\end{lemma}
\proof
Let us first recall the so called Jonqui\`ere's inversion formula
\cite{Jo} (see also Appendix \ref{app:1}), which provides a description
of the analytic continuation of the polylogorithm functions in terms
of Bernoulli polynomials
\beq\label{poly-rec}
\operatorname{Li}_p(1/x) = (-1)^{p+1} \operatorname{Li}_p(x)+
(-1)^{p+1} \frac{(2\pi\sqrt{-1})^p}{p!} \,
B_p\Big(\frac{1}{2\pi\sqrt{-1}}\log x\Big),\quad 0<|x| < 1,
\eeq
where $B_p(x)$ $(p\geq 0)$ are the Bernoulli polynomials defined by 
\beq\label{Bernoulli}
\frac{t e^{xt}}{e^t-1} = \sum_{p=0}^\infty B_p(x)\frac{t^p}{p!}.
\eeq
The value of $\operatorname{Li}_p(1/x)$ is specified via analytic
continuation along a path $C'\subset \mathbb{C}-\{1\}$ from $x$ to
$1/x$, which {\em does not} wind
around $x=1$.  The choice of the branch of $\log x$ in
\eqref{poly-rec} is such that the formula holds for $p=1$. Let us
assume first that $C'$ intersects the real axis once and that
$\operatorname{Im}(x)\neq 0$. In order to
determine the branch of $\log x$, we have to consider 4 cases
depending on whether $\operatorname{Im}(x)>0$ or $<0$ and whether $C'$
intersects the real interval $(1,\infty)$ or not. However, after
analyzing the 4 cases we find that there are two possibilities only. Namely, if 
we walk along $C'$ from $x$ to $1/x$, then when crossing the real
axis, 1 is either on our left or on our right. 
In the first case $\log x:=\operatorname{Log} x $, otherwise  $\log
x:=\operatorname{Log} x +2\pi\sqrt{-1} $, where  
\ben
\operatorname{Log} x := \ln |x|+\sqrt{-1}\, \operatorname{Arg}(x),\quad 
-\pi<\operatorname{Arg}(x)\leq \pi
\een
is the principal branch of the logarithm. For general $C'$, the
choice of the branch can be deduced  easily from the above two cases.  

Put $x=\mu/\lambda$. 
Recalling formula \eqref{poly} and using \eqref{poly-rec} we get
\ben
\operatorname{Li}_\sigma(1/x) = 
\operatorname{Li}_\sigma(x)^T+
\sum_{r=1}^{|\sigma|}\sum_{p=1}^\infty
\frac{x^{-\mathcal{N}_n}}{\mathcal{N}_n|\sigma|}\,
\frac{(2\pi\sqrt{-1} \mathcal{N}_n|\sigma|)^p}{p!}\,
B_p\Big(\frac{1}{2\pi\sqrt{-1}}\log x_r\Big)\sigma_s^{-r},
\een
where $x_r=\eta^r x^{1/|\sigma|}$ and the paths $C_r'$ from $x_r$ to
$x_r^{-1}$ can be described as follows. Let $C': t\mapsto
C_2(t)/C_1(t)$ be the path from $x$ to $1/x$ and $C'_0$ be the induced
path from $x^{1/|\sigma|}:=e^{\operatorname{Log}(x)/|\sigma|}$ to
$x^{-1/|\sigma|}:=e^{-\operatorname{Log}(x)|/\sigma|}$. Then $C'_r$ is
a composition of a path inside the unit disk from $x_r$ to $x_{-r}$
and $\eta^{-r}C'_0$ (clockwise rotation of $C'_0$ with angle $2\pi
r/|\sigma|$).  

The infinite sum over $p$ can be computed via \eqref{Bernoulli}
\ben
\operatorname{Li}_\sigma(1/x) = 
\operatorname{Li}_\sigma(x)^T+
\sum_{r=1}^{|\sigma|}
\frac{x^{-\mathcal{N}_n}}{\mathcal{N}_n|\sigma|}\,\Big(
2\pi\sqrt{-1} \mathcal{N}_n|\sigma|\frac{e^{\mathcal{N}_n|\sigma \log x_r} } {e^{2\pi\sqrt{-1}
    \mathcal{N}_n |\sigma|} - 1 } -1\Big)\sigma_s^{-r}.
\een
Since
$(0,\lambda,\mu)\in D_\infty^+$ we have $|x-1|<c<1$ and
$|x|<1$. Moreover, since the branch of the period vectors depending on
$\mu$  is determined from the branch of the period vectors depending
on $\lambda$ via the straight segment from $[\lambda,\mu]$, we get
that the multi-valued functions
\ben
x^{-\mathcal{N}_n } = e^{-\mathcal{N}_n\operatorname{Log}x },\quad 
x^{1/|\sigma|} = e^{\frac{1}{|\sigma|} \, \operatorname{Log}x }
\een
are defined via the principal branch of the logarithm. 
Hence 
\ben
x^{-\mathcal{N}_n} e^{\mathcal{N}_n|\sigma| \log x_r} = 
\exp\Big( \mathcal{N}_n |\sigma|\sqrt{-1} \Big( 
\operatorname{Arg}(\eta^rx^{1/|\sigma|})-
\operatorname{Arg}(x^{1/|\sigma|})+
2\pi\chi_r \Big)\Big),
\een
where $\chi_r=1$ or $0$, depending whether $1$ is on the left or
on the right of $C'_r$.
Since $C_0'$ is in a small neighborhood of $1$ we have
\ben
\operatorname{Arg}(\eta^rx^{1/|\sigma|})-
\operatorname{Arg}(x^{1/|\sigma|})+
2\pi\chi_r  = 
2\pi r/|\sigma|,\quad 1\leq r\leq |\sigma|-1.
\een
The sum over $r$ becomes 
\ben
\frac{2\pi\sqrt{-1}}{\sigma^{-|\sigma|}-1}\,\Big(
\sum_{r=1}^{|\sigma|-1}
\sigma^{-r} + \sigma^{-|\sigma|\chi_0}\Big) +
\frac{x^{-\mathcal{N}_n}}{-\mathcal{N}_n}  \Big(\frac{1}{|\sigma|}
\sum_{r=1}^{|\sigma|}\sigma_s^r\Big).
\een
Finally, we get
\beq\label{Lis-rec}
\operatorname{Li}_\sigma(1/x) = 
\operatorname{Li}_\sigma(x)^T-2\pi\sqrt{-1}\, \frac{\sigma^{1-\chi_0}
}{\sigma-1}+
\frac{x^{-\mathcal{N}_n}}{-\mathcal{N}_n}  \Big(\frac{1}{|\sigma|}
\sum_{r=1}^{|\sigma|}\sigma_s^r\Big). 
\eeq
Using the formula for the phase factors in Theorem \ref{t1} we compute
\ben
\Omega_{\alpha,\beta}(0,\lambda,\mu)-\Omega_{\beta,\alpha}(0,\mu,\lambda)
=\Big( ( \operatorname{Li}_\sigma(1/x) -\operatorname{Li}_\sigma(x)^T)
\beta|\alpha\Big) +P_{\alpha,\beta}(\lambda,\mu)-P_{\beta,\alpha}(\mu,\lambda).
\een
Recalling again Theorem \ref{t1} and using \eqref{Lis-rec} we get 
\ben
\Omega_{\alpha,\beta}(0,\lambda,\mu)-\Omega_{\beta,\alpha}(0,\mu,\lambda)
= 2\pi\sqrt{-1}\operatorname{SF}(\sigma^{1-\chi_0}\beta,\alpha).
\een
If $\chi_0=0$, then
$\operatorname{SF}(\sigma\beta,\alpha)=-\operatorname{SF}(\alpha,\beta)$
otherwise, if $\chi_0=1$, then 
$\operatorname{SF}(\beta,\alpha)=-\operatorname{SF}(\alpha,\beta)+(\alpha|\beta)$. 
\qed

\subsection{Analytic extension of the phase factors}

Recall that the series \eqref{inf-series} is convergent in the domain
\ben
D^+:=
\{(t,\lambda,\mu)\ |\ |\lambda-\mu|<|\mu|-r(t)<|\lambda|-r(t)\}
\een
It is not hard
to see that the RHS of \eqref{pf-partial} has a singularity at
$\lambda=\mu$ of the form $(\alpha|\beta)(\lambda-\mu)^{-1}$. Hence
the Laurent series expansion at $\lambda=\infty$ of the series 
\ben
\widetilde{\Omega}_{\alpha,\beta}(t,\lambda,\mu) :=
 {\Omega}_{\alpha,\beta}(t,\lambda,\mu) - (\alpha|\beta)\log (\lambda-\mu)
\een
is convergent in the domain
\ben
D=\{(t,\lambda,\mu)\ |\ |\lambda-\mu|<
\operatorname{min}\{|\mu|-r(t),|\lambda|-r(t)\}\}.
\een
In particular, we get that ${\Omega}_{\alpha,\beta}(t,\lambda,\mu)$
extends analytically along any path in the domain $D_\infty-\{\lambda=\mu\}$ and that
$e^{{\Omega}_{\alpha,\beta}(t,\lambda,\mu)}$ has a pole of order
$-(\alpha|\beta)$ at $\lambda=\mu$. 
\begin{lemma}\label{convergence}
If $(t,\lambda,\mu)\in D_\infty^+$ and $C\subset D_\infty-\{\lambda=\mu\}$ is a 
path from $(t,\lambda,\mu)$ to $(t,\mu,\lambda)$, s.t., the projections of
$C$ along the last two coordinates are homotopic to the 
line segments $[(t,\lambda),(t,\mu)]$ and $[(t,\mu),(t,\lambda)]$, then  
\ben
\Omega_{\alpha,\beta}(t,\lambda,\mu)-\Omega_{\beta,\alpha}(t,\mu,\lambda)=
- 2\pi\sqrt{-1} \Big(
\operatorname{SF}(\alpha,\beta)+k(\alpha|\beta)
\Big),\quad k\in \mathbb{Z},
\een
where the 2nd phase factor on the LHS is obtained from
$\Omega_{\beta,\alpha}(t,\lambda,\mu)$ via analytic continuation along
the path $C$.
\end{lemma}
\proof
The main idea is to obtain an integral representation of the infinite series 
$\Omega_{\alpha,\beta}(t,\lambda,\mu)$. 
If $(t,\lambda,\mu)\in D_\infty$, then by definition we can find a path $C\subset
S_\infty$ from $(0,\lambda)$ to $(t,\lambda)$, s.t.,
$|\mu-\lambda|<|\lambda'|-r(t')$ for all $(t',\lambda')\in C$. Put
$\mu':=\lambda'+\mu-\lambda$. Using the triangle inequality we get
$|\mu'|\geq |\lambda'|-|\mu-\lambda|>r(t')$. Therefore,
$(t',\lambda',\mu')\in D$ and it makes sense to think of
$\widetilde{\Omega}_{\alpha,\beta}(t',\lambda',\mu')$ as a holomorphic function on
the curve $C$.  Using the differential equations
\eqref{PF-1}--\eqref{de-prim} we get 
\ben
d \widetilde{\Omega}_{\alpha,\beta}(t',\lambda',\mu') =
\left. \tau^*\mathcal{W}_{\alpha,\beta}\right|_C,
\een
where $\tau:S\to B$, $(t,\lambda)\mapsto
t-\lambda\mathbf{1}$. Integrating the above identity along $C$ we get  
\beq\label{phase-int}
\Omega_{\alpha,\beta}(t,\lambda,\mu)
=\Omega_{\alpha,\beta}(0,\lambda,\mu) +
\int_{-\lambda\mathbf{1}}^{t-\lambda\mathbf{1}}
\mathcal{W}_{\alpha,\beta}(t'-\lambda'\mathbf{1}, \mu-\lambda).
\eeq
The translation $C+(0,\mu-\lambda)$ is a path $\widetilde{C}$ in $S_\infty$
from $(0,\mu)$ to $(t,\mu)$. Note that $\widetilde{C}$ is obtained
from $C$ also by the map $(t',\lambda')\mapsto
(t',\mu')$. 
Using that by definition 
\ben
\mathcal{W}_{\alpha,\beta}(t'-\lambda'\,\mathbf{1},\mu-\lambda) =
I^{(0)}_\alpha(t',\lambda')\bullet_{t'} I^{(0)}_\beta(t',\mu'),
\een
and translating the integration variables $(t',\lambda')\mapsto
(t',\mu')$ we get
\ben
\int_{-\lambda\mathbf{1}}^{t-\lambda\mathbf{1}}
\mathcal{W}_{\alpha,\beta}(t'-\lambda'\mathbf{1}, \mu-\lambda)=
\int_{-\mu\mathbf{1}}^{t-\mu\mathbf{1}}
\mathcal{W}_{\beta,\alpha}(t'-\mu'\mathbf{1}, \lambda-\mu),
\een
i.e., the integral in \eqref{phase-int}
is symmetric with respect to exchanging the pairs
$(\alpha,\lambda)$ and $(\beta,\mu)$. To finish the proof it
remains only to recall Lemma \ref{phase-sym:1}.\qed

\subsection{Proof of Theorem \ref{t2}}
Let us choose $t\in B$, s.t., the canonical coordinates $u_i:=u_i(t)$
are vertices of a regular $N$-gon with center at $0$. In particular,
$r(t)=|u_i|$ for all $i$. Let us fix a
reference point $(t,\lambda)$ in $\{t\}\times \mathbb{C}\subset S_\infty$ and take a closed loop
$C$ obtained by approaching one of the points $u_i$ along a path in
$S_\infty$ and going around $u_i$. Since the fundamental group of $S'$
is generated by loops of the above type, it suffices to prove that 
\beq\label{ident}
\oint_{ C}
\tau^*\mathcal{W}_{\alpha,\beta}(\, \cdot\, ,\mu-\lambda) -
\Omega_{w(\alpha),w(\beta)}(t,\lambda,\mu)+
\Omega_{\alpha,\beta}(t,\lambda,\mu) 
\quad \in \quad 
2\pi\sqrt{-1}\mathbb{Z},
\eeq
where $w$ is the monodromy
transformation along $C$ and $\mu$ is sufficiently close to
$\lambda$. The integral in \eqref{ident} can be written also as
follows: if $(t,x)\in C$, then $\tau(t,x)=t-x\mathbf{1} \in \tau(C)$
and we have
\ben
\oint_C \tau^*\mathcal{W}_{\alpha,\beta}(\, \cdot\, ,\mu-\lambda)=
\oint_{\tau(C)}\mathcal{W}_{\alpha,\beta}(t-x\mathbf{1},\mu-\lambda).
\een
In our argument we always keep $t$ fixed,
all paths in $S$ that we use will be identified with paths in
$\mathbb{C}\cong \{t\}\times \mathbb{C}\subset S$, and we drop
the argument of the phase form. 

\begin{lemma}\label{vv}
If $\varphi$ is the cycle vanishing over $u_i$, then $\oint_{ C}
\tau^*\mathcal{W}_{\varphi,\varphi}= -4\pi\sqrt{-1}$. 
\end{lemma}  
\proof
The differential equations of the periods imply that the integral is
independent of $\lambda$ and $\mu$. Therefore, we may set
$\mu=\lambda$. The computation in this case can be found in
\cite{G3}. 
\qed

\begin{lemma}\label{iv}
Assume that $\varphi$ is the cycle vanishing over $u_i$ and $\alpha$
is a cycle invariant along $C$, then 
\ben
\oint_{ C}
\tau^*\mathcal{W}_{\alpha,\varphi}=-2\Omega_{\alpha,\varphi}(t,\lambda,\mu).
\een
\end{lemma}
\proof
Let us choose $\lambda'\in C$ sufficiently close to $u_i$ with $|\lambda'|>r(t)$ and put
$\mu'=\lambda'+\mu-\lambda$. The infinite series \eqref{inf-series}
at $(t,\lambda,\mu)=(t,\lambda',\mu')$ can be interpreted in two
different ways: as a convergent Laurent series in $(\lambda')^{-1}$ or as a
convergent Laurent series in $(\mu'-u_i)^{1/2}$. The 2nd
interpretation allows us to extend the equality 
$d\Omega_{\alpha,\varphi}(t,\lambda',\mu')=\tau^*\mathcal{W}_{\alpha,\varphi}$
as $\lambda'$ varies along the loop around $u_i$. The Lemma follows by
the Stoke's theorem.\qed

If $\alpha$ and $\beta$ are arbitrary cycles, then we can decompose
them as
\ben
\alpha=\alpha'+(\alpha|\varphi)\varphi/2,\quad 
\beta=\beta'+(\beta|\varphi)\varphi/2,
\een
where $\alpha'$ and $\beta'$ are invariant along $C$. We have to
compute
\ben
\oint_C \tau^* \mathcal{W}_{\alpha,\beta} = 
\frac{(\beta|\varphi)}{2}
\oint_C \tau^* \mathcal{W}_{\alpha',\varphi}+
\frac{(\alpha|\varphi)}{2}
\oint_C \tau^* \mathcal{W}_{\varphi,\beta'}+
\frac{(\alpha|\varphi) (\beta|\varphi)}{4}
\oint_C \tau^* \mathcal{W}_{\varphi,\varphi}
\een 
where we used that $\oint_C \tau^*\mathcal{W}_{\alpha',\beta'}=0$. The
last and the first integrals can be computed by respectively Lemma
\ref{vv} and Lemma \ref{iv}. For the middle one, note that changing
the variables $(t,x)\to (t,x+\mu-\lambda)$ gives an integral along a
loop $\widetilde{C}$ based at $(t,\mu)$ which can be computed again by Lemma \ref{iv}
\ben
\oint_C \tau^* \mathcal{W}_{\varphi,\beta'}(\cdot,\mu-\lambda) = 
\oint_{\widetilde{C}} \tau^*
\mathcal{W}_{\beta',\varphi}(\cdot,\lambda-\mu) = 
-2\Omega_{\beta',\varphi}(t,\mu,\lambda). 
\een
In other words, we have
\ben
\oint_C \tau^* \mathcal{W}_{\alpha,\beta} = 
-(\beta|\varphi)\Omega_{\alpha',\varphi}(t,\lambda,\mu)
-(\alpha|\varphi)\Omega_{\beta',\varphi}(t,\mu,\lambda)
- \pi\sqrt{-1} (\alpha|\varphi) (\beta|\varphi).
\een
Therefore, \eqref{ident} takes the form
\ben
(\alpha|\varphi)
\Big(
\Omega_{\varphi,\beta'}(t,\lambda,\mu) - 
\Omega_{\beta',\varphi}(t,\mu,\lambda)
\Big)
- \pi\sqrt{-1} (\alpha|\varphi) (\beta|\varphi).
\een
Recalling Lemma \ref{convergence} we get 
\ben
&&
-2\pi\sqrt{-1} (\alpha|\varphi) \operatorname{SF}(\varphi,\beta')
- \pi\sqrt{-1} (\alpha|\varphi) (\beta|\varphi)=\\
&&
=
-2\pi\sqrt{-1} (\alpha|\varphi) \operatorname{SF}(\varphi,\beta) 
+\pi\sqrt{-1} (\alpha|\varphi) (\beta|\varphi)\Big(
\operatorname{SF}(\varphi,\varphi) - 1\Big)
.
\een
Since $\operatorname{SF}(\varphi,\varphi)=1$, the above number is an
integer multiple of $2\pi\sqrt{-1}$. \qed

\section{VOA representations}\label{sec:VOA}
In this section we would like to explain the origin of the phase
factors and the importance of Theorem
\ref{t2} from the point of view of the representation theory of lattice VOAs.

\subsection{The tame Fock space}

Let us denote by 
\ben
\mathfrak{W}:=\Big(\bigoplus_{m=0}^\infty \mathcal{H}^{\otimes m} \otimes
\mathbb{C}(\!(\hbar )\!)\Big) / I
\een
the Weyl algebra of $\mathcal{H}$, where $I$ is the two-sided ideal
generated by 
\ben
v\otimes w-w\otimes v -\Omega(v,w)\hbar,\quad v,w\in
\mathcal{H}.
\een 
Furthermore we introduce the {\em Fock space} 
$\mathbb{V} = \mathfrak{W}/\mathfrak{W} \mathcal{H}_+$.
For any $\mathbb{C}$-algebra $A$, let us denote by $\widehat{\mathbb{V}}_A$ the completion of
$\mathbb{V}\otimes A$ corresponding to the filtration $\{\mathcal{H}_-^m
\mathbb{V}\otimes A\}_{m=0}^\infty$. When $A=\mathbb{C}$ then the
index $A$ will be omitted.

As a vector space the Weyl algebra $\mathfrak{W} $ is isomorphic to
\ben
\bigoplus_{m',m''=0}^\infty \mathcal{H}_-^{\otimes
  m'}\otimes \mathcal{H}_+^{\otimes m''}\otimes\,  \mathbb{C}(\!(\hbar)\!).
\een
We introduce the {\em tame Weyl algebra}\footnote{This notion was
  invented by H. Iritani}
\ben
\mathfrak{W}_{\rm tame} \subset \widehat{\mathfrak{W}}:=\prod_{m',m''=0}^\infty
\Big(
\underbrace{\mathcal{H}_-\widehat{\otimes} \cdots \widehat{\otimes} \mathcal{H}_- }_{m'}\widehat{\otimes}
\underbrace{\mathcal{H}_+\widehat{\otimes} \cdots \widehat{\otimes} \mathcal{H}_+}_{m''}
\widehat{\otimes} \, \mathbb{C}(\!(\hbar)\!)
\Big),
\een
where $\widehat{\otimes}$ is the completion of the tensor product
induced by the filtrations $\{z^{-m}\mathcal{H}_-\}_{m=0}^\infty$,
$\{z^m\mathcal{H}_+\}_{m=0}^\infty$, and
$\{\hbar^m\mathbb{C}(\!(\hbar)\!)\}_{m=0}^\infty$ of respectively
$\mathcal{H}_-$, $\mathcal{H}_+$, and $\mathbb{C}(\!(\hbar)\!)$. Every
element of $\widehat{\mathfrak{W}}$ is an infinite series of monomials
of the type $\hbar^{g-1} v\otimes w$,
where $g\in \mathbb{Z}$, $v$ and $w$ have the form respectively $v_1 z^{-k_1-1}\otimes \cdots \otimes
v_{m'}z^{-k_{m'}-1}$ and $w=w_1z^{\ell_1}\otimes\cdots\otimes
w_{m''}z^{\ell_{m''}}$. The tame Weyl algebra $\mathfrak{W}_{\rm
  tame}$ consists of those series for which the non-zero monomials
satisfy the inequality
\ben
k_1+\cdots+k_{m'} - m'\leq 3 (g-1 + m''/2).
\een
Let us point out that if $m''=0$, then this is the notion of a tame
function introduced in \cite{G3}. One can check that $\mathfrak{W}_{\rm
  tame}$ is in fact an algebra with the product induced from
$\mathfrak{W}$. 
Finally, let us introduce also the {\em tame } Fock
space $\mathbb{V}_{\rm tame}:=\mathfrak{W}_{\rm
  tame}/\mathfrak{W}_{\rm tame}\mathcal{H}_+$.

\subsection{Heisenberg fields and Vertex operators}

Let us denote by $\mathcal{O}$ the algebra of analytic functions on
the monodromy covering space $\widetilde{S'}$ with at most polynomial
growth at $\lambda=\infty$. By definition (and since the
Gauss--Manin connection has regular singularities) all
period vectors $I^{(m)}_\alpha\in \mathcal{O}\otimes H$. Let us
define the following set of linear operators $\mathbb{V}_{\rm tame}\to
\widehat{\mathbb{V}}_\mathcal{O}$:
\beq\label{bosonic-field}
\phi_a(t,\lambda) = \hbar^{-1/2} \partial_\lambda
\mathbf{f}_a(t,\lambda;z),\quad a\in \mathfrak{h}^*
\eeq
and 
\beq\label{vop}
\Gamma^\alpha (t,\lambda) =
e^{\hbar^{-1/2}\mathbf{f}_\alpha(t,\lambda;z)_-}e^{\hbar^{-1/2}\mathbf{f}_\alpha(t,\lambda;z)_+},
\quad    
\alpha\in Q,
\eeq
where $\pm$ denotes the projection on $\mathcal{H}_\pm$. If we want to
define a VOA-representation using the above operators we need to
define a composition of vertex operators. Let us first compose two
vertex operators ignoring the possible divergence issues. We get 
\beq\label{voa-prod}
\Gamma^\alpha(t,\lambda)\Gamma^\beta(t,\mu) =
e^{\Omega_{\alpha,\beta}(t,\lambda,\mu)}
:\Gamma^\alpha(t,\lambda)\Gamma^\beta(t,\mu):,
\eeq
where $:\ :$ is the normally ordered product for which the action of all elements of
$\mathcal{H}_+$ precedes the action of $\mathcal{H}_-$. 
The problem is that in general we can not exponentiate elements of
the ring
$\mathbb{C}(\!(\mu^{-1/|\sigma|})\!)(\!(\lambda^{-1/|\sigma|})\!) $. 
To offset this difficulty we recall that the series
\eqref{inf-series} is convergent for all $(t,\lambda,\mu)\in
D_\infty-\{\lambda=\mu\}$. Moreover, using the integral presentation
\eqref{phase-ac} we can extend analytically
$\Omega_{\alpha,\beta}(t,\lambda,\mu)$ with respect to $(t,\mu)\in S'$
provided $\lambda$ is sufficiently close to $\mu$. 
Therefore, the {\em phase factor}
\beq\label{phase-factor}
B_{\alpha,\beta}(t,\lambda,\mu):=\exp\left( \Omega_{\alpha,\beta}(t,\lambda,\mu) \right)
\eeq
has a convergent Laurent series
expansion in $(\mu-\lambda)$  with coefficients multi-valued analytic
functions on $S'$, i.e., the coefficients are analytic on
the universal cover of $S'$. According to Theorem \ref{t2} 
the coefficients are analytic on the {\em monodromy} covering space
$\widetilde{S'}$.  In other words we have the following proposition. 
\begin{proposition}\label{phase-factor-ac}
The phase factor $B_{\alpha,\beta}(t,\lambda,\mu)$ has a Laurent series
expansion in $\mu-\lambda$, whose coefficients are elements of
$\mathcal{O}$.  
\end{proposition}

\subsection{Applications to VOA}

Let us denote by 
\ben
V_Q:=\operatorname{Sym}(\mathfrak{h}^*[s^{-1}]s^{-1})\otimes \mathbb{C}[Q]
\een 
the vector space underlying the lattice vertex algebra corresponding to the
Milnor lattice $Q$. Note that we did
not put a VOA structure on $V_Q$ yet. Following \cite{BM}, we can
introduce a family of products $a_{(k)}b$, $a,b\in V_Q$, $k\in
\mathbb{Z}$ and a linear map
\ben
X: V_Q\to \operatorname{Hom}_\mathbb{C} (\mathbb{V}_{\rm
  tame},\widehat{\mathbb{V}}_\mathcal{O}), 
\een
defined recursively by the initial conditions
\beq\label{initial}
X(as^{-1}\otimes 1) := \phi_a,\quad X(1\otimes e^\alpha) := \Gamma^\alpha,
\eeq
and the {\em operator product expansion} (OPE) formula
\beq\label{ope}
X_t(a_{(M-k-1)}b,\lambda) \mathcal{A} = \frac{1}{k!}\partial_\mu^k\left.\left( 
(\mu-\lambda)^M X_t(a,\mu)X_t(b,\lambda)\mathcal{A}\right)\right|_{\mu=\lambda},
\eeq
where $a,b\in V_Q$ and $\mathcal{A}\in \mathbb{V}_{\rm tame}$ are
arbitrary, the number $M$ is sufficiently large, $k\geq 0$, and we
denote by $X_t(a,\lambda)$ the value of the operator 
$X(a)$ at a point $(t,\lambda)\in S'$. One has to check that the above
definition is correct and that all the fields $X_t(a,\lambda)$, $a\in
V_Q$, are mutually local. Let us point out that the products
$a_{(k)}b$ will depend on the choice of a point on $\widetilde{S'}$,
i.e., a point $(t,\lambda)$ with a reference path. Our expectation is
that $V_Q$ has a unique structure of a VOA, s.t., the products defined
via the OPE formula coincides with the Borcherds' products. Moreover,
the VOA structure on $V_Q$ is isomorphic to the lattice VOA structure
(see \cite{Kac} for some background). The details will be presented
elsewhere.

\appendix
\section{Jonqui\`ere's inversion formula}\label{app:1}

We would like to give a proof of formula \eqref{poly-rec}, because as
it was pointed out by the Wikipedia article about polylogorithms, some
authors are not very careful about the choice of the branch of $\log
x$, while for our purposes choosing the correct branch is essential. 

Let us assume that $C'$ intersects the real axis once and that ${\rm
  Im}(x)\neq 0$. We can always achieve this by replacing $C'$ with a
homotopy equivalent path and by deforming slightly $x$. Put 
\ben
f_p(x) := \operatorname{Li}_p(x)+(-1)^p  \operatorname{Li}_p(1/x).
\een  
Since $x\partial_x f_p(x) = f_{p-1}(x)$, we get that $f_p(x)$ is a polynomial in
$\log x$. Let us write it as 
\ben
f_p(x)=:-\frac{(2\pi\sqrt{-1})^p}{p!}h_p\Big( \frac{\log
  x}{2\pi\sqrt{-1}}\Big).
\een 
In order to determine the polynomials $h_p$, let us see what happens when we
analytically continue $x$ along a closed loop inside the unit disk going
counterclockwise around $0$.  There are 4 cases to be considered,
depending on whether $C'$ intersects
$(1,\infty)$ or not and whether ${\rm Im}(x)>0$ or $<0$. Let us
analyze the case: $C'$ does not intersect $(1,\infty)$ and ${\rm
  Im}(x)<0$. The remaining cases are similar. Note that in this case
$\log x = \operatorname{Log} x +2\pi\sqrt{-1}$. We have the following
integral representation
\ben
\operatorname{Li}_p(x^{-1}) = \frac{x^{-1}}{\Gamma(p)} 
\int_1^\infty \frac{(\operatorname{Log}\, z)^{p-1}}{z-x^{-1}}\frac{dz}{z}.
\een
When $x$ is varying around $0$, the point $x^{-1}$ will vary clockwise along a
closed loop $C$ around the unit disk. Let $y$ be a point on $C$ which
we reach slightly after we cross the interval $(1,\infty)$. Using the
Cauchy residue theorem,  we get
\ben
\operatorname{Li}_p(y) = \frac{y}{\Gamma(p)} 
\int_1^\infty \frac{(\operatorname{Log}\, z)^{p-1}}{z-y}\frac{dz}{z} + 
\frac{2\pi\sqrt{-1}}{\Gamma(p)} 
(\operatorname{Log}\, y)^{p-1}.
\een
Extending analytically the above identity as $y$ varies along the
remaining part of $C$, we get the following formula
\ben
{\rm a.c.} f_p(x)-f_p(x)=
(-1)^p\Big({\rm a.c.} \operatorname{Li}_p(x^{-1})  -
\operatorname{Li}_p(x^{-1}) \Big) = 
-\frac{2\pi\sqrt{-1}}{\Gamma(p)} 
\Big(-\operatorname{Log}\, x^{-1} +2\pi\sqrt{-1} \Big)^{p-1}.
\een
Expressing $f_p$ in terms of $h_p$ we arrive
at the following difference equation
\ben
h_p(L+1)-h_p(L) = p L^{p-1},\quad L:=\log x /2\pi\sqrt{-1},\quad p\geq 1.
\een
On the other hand, since $x\partial_x f_p(x)=f_{p-1}(x)$ we also have
$\partial_Lh_p(L)=ph_{p-1}(L)$. The solution to these relations is
unique and it is given by the Bernoulli polynomials.

\section*{Acknowledgements} 
I am thankful to K. Hori, C. Li, Si Li, and K. Saito for organizing
the workshop ``Primitive Forms and Related Subjects'' (Kavli IPMU,
Feb/2014) and creating an
inspiring environment. I am especially thankful to S. Galkin and
H. Iritani for pointing out that Dubrovin's classification of
semi-simple Frobenius manifolds can be used to extend the domain of
the period vectors. This work is supported by JSPS Grant-In-Aid 26800003 and by the World
Premier International Research Center Initiative (WPI Initiative), MEXT, Japan.

\end{document}